\documentclass[reqno]{amsart}

\usepackage{graphicx}
\usepackage{times}
\usepackage{amsmath}
\usepackage{amsthm}
\usepackage{amssymb,color}
\usepackage{lineno,hyperref}
\usepackage{fancyhdr}
\usepackage{fancyheadings}
\usepackage{paralist}
\usepackage{bbm} 

\newtheorem{theorem}{Theorem}[section]
\newtheorem{atheorem}{Theorem}

\theoremstyle{definition}
\newtheorem{definition}[theorem]{Definition}

\theoremstyle{remark}

\numberwithin{equation}{section}

\begin{document}

\title[Stable and Historic Behavior in Replicator Equations]{Stable and Historic Behavior in Replicator Equations Generated by Similar-Order Preserving Mappings}

\author{Mansoor Saburov}
\address{Mansoor Saburov: College of Engineering and Technology, American University of the Middle East, Kuwait}
\email{mansur.saburov@aum.edu.kw}


\date{}


\begin{abstract}
	One could observe drastically different dynamics of \textit{zero-sum} and \textit{non-zero-sum} games under \textit{replicator equations}. In \textit{zero-sum} games, \textit{heteroclinic cycles} naturally occur whenever the spe\-cies of the population supersede each other in cyclic fashion (like for the \textit{Rock-Paper-Scissors} game). In this case, the highly erratic oscillations may cause \textit{the divergence of the time averages}. In contrast, it is a common belief that the most \textit{``reasonable"} replicator equations of \textit{non-zero-sum} games satisfy \textit{``The Folk Theorem of Evolutionary Game Theory"} which asserts that $(i)$ \textit{a Nash equilibrium is a rest point}; $(ii)$ \textit{a stable rest point is a Nash equilibrium}; $(iii)$ \textit{a strictly Nash equilibrium is asymptotically stable}; $(iv)$ \textit{any interior convergent orbit evolves to a Nash equilibrium}.  In this paper, we propose two distinct vast classes of replicator equations generated by \textit{similar-order preserving mappings} which exhibit \textit{stable} as well as \textit{mean historic behavior}. In the latter case, the time averages of the orbit will slowly oscillate during the evolution of the system and do not converge to any limit. This will eventually cause the divergence of higher-order repeated time averages.
\end{abstract}

\maketitle



\section{Introduction}

The ergodic theory studies 
the long term average 
behavior of systems evolving in time.
Two of the most important theorems are
those of Birkhoff and von Neumann 
which assert the convergence of time averages
along almost each orbit.
For the special class of ergodic systems,
these time averages have the same limit 
for almost all initial points:
statistically speaking, 
\textit{the system that evolves 
	for a long time forgets its initial state},
i.e., \textit{it does not have  history}.
The set of all initial points 
for which the time averages 
do not converge has measure 
zero with respect to 
any invariant measure.
This however does not imply
that it has Lebesgue measure zero, 
in general. 
Furthermore, even if the set has 
Lebesgue measure zero,
it is possible for the set 
to be topologically large \cite{Dowker,Thompson2010,Thompson2012}.
This motivates to introduce a notion
of \textit{historic behavior} 
in a dynamical system.
Roughly speaking, we say that 
the forward orbit has 
\textit{historic behavior} 
if its time averages do not converge 
in the embedded topology. The terminology 
\textit{historic behavior} was coined
by Ruelle \cite{Ruelle}.
The problem of describing 
the persistent family of 
dynamical systems with historic behavior 
was popularized by Takens 
\cite{Takens2008}. 
Recently, this problem 
under the name of 
\textit{Takens' Last Problem}
was studied 
in the papers 
\cite{BKNRS,MCPV2021,deSantana2022,KLS16,KNS19,KNS22,KS17,Yang2020}.

These two extreme behaviors, \textit{stable} and \textit{historic behavior}, of dynamical systems can be frequently observed in evolutionary game theory. The main aim of this paper is to provide two vast classes of replicator equations which exhibit \textit{stable} as well as \textit{historic behavior}. 

Evolutionary game theory  is the study of \textit{frequency-dependent natural selection} in which the fitness of individuals is \textit{not constant}, but depends on \textit{frequencies} of the different phenotypes in population \cite{MaynardSmithPrice}. Evolutionary game theory differs from classical non-cooperative game theory in focusing more on \textit{the dynamics of strategy change}. Moreover, it does not rely on rationality, i.e., the players do not choose their strategy and cannot change it: they are born with a strategy and their offspring inherit that same strategy. Evolutionary game theory interprets payoff as \textit{biological fitness} or \textit{reproductive success}. The success of a strategy is determined by how good the strategy is in the presence of competing strategies. The strategies that do well reproduce faster and  the strategies that do poorly are out-competed. There has been a veritable explosion of interest by economists and social scientists in evolutionary game theory \cite{Friedman1998,PST2022,Sandholm,Sigmund2010,SunHouSun2022}. 

\textit{The Nash equilibrium} is a strategy profile in which no player can do better by unilaterally changing their strategy to another strategy \cite{Nash1951}. An \textit{evolutionarily stable strategy} (ESS) is akin to the Nash equilibrium which is \textit{``evolutionarily" stable}, i.e., once it is fixed in a population, natural selection alone is sufficient to prevent \textit{mutant} strategies from invading successfully \cite{Thomas1984,Thomas1985}. Therefore, the ESS is effective against any \textit{mutant} strategy when it is initially rare and successful when it is eventually abundant.

\textit{The replicator equation} is the cornerstone of evolutionary game theory  \cite{TaylorJonker}. A \textit{replicator} can be a \textit{gene} in population genetics, a \textit{molecule} in biochemical evolution, an \textit{organism} in population ecology, a \textit{strategy} in evolutionary game. Several evolutionary models in distinct biological fields -- \textit{population genetics},  \textit{biochemical evolution}, \textit{population ecology}, and \textit{evolutionary game} -- lead independently to the same class of replicator equations \cite{SchusterSigmund,TaylorJonker}. The general idea is that replicators whose fitness is larger (smaller) than the average fitness of population will increase (decrease) in numbers. In replicator equation,  the relative growth rate of players using a certain strategy is equal to the difference between the average payoff of that strategy and the average payoff of the population as a whole. Hence, the static approach of evolutionary game theory has been complemented by the dynamic stability analysis of \textit{rest} (\textit{stationary}) solutions of the replicator equations.

On the one hand, it is a common belief \cite{Cressman2003,HofbauerSigmund1998,HofbauerSigmund2003} that the most \textit{``reasonable"} replicator equations satisfy the \textit{``Folk Theorem of Evolutionary Game Theory"}  which asserts that $(i)$ \textit{a Nash equilibrium is a rest point}; $(ii)$ \textit{a stable rest point is a Nash equilibrium}; $(iii)$ \textit{a strictly Nash equilibrium is asymptotically stable}; $(iv)$ \textit{any interior convergent orbit evolves to a Nash equilibrium}. 

On the other hand, in the zero-sum game \cite{Gaunersdorfer1992,Hofbauer,Sigmund}, heteroclinic cycles naturally occur whenever the spe\-cies of the population supersede each other in cyclic fashion (like for the \textit{Rock-Paper-Scissors} game). The dynamics near the stable heteroclinic cycle is characterized by \textit{intermittency}: an orbit remaining near the heteroclinic cycle, spends long periods of time close to each saddle and makes fast transitions from one saddle to the next one. Since the orbit converges to the heteroclinic cycle, it comes closer and closer to the saddles and consequently remains there for longer and longer times and the time spent near each saddle increases \textit{exponentially}. Such highly erratic oscillations cause \textit{the divergence of the time averages} which is known as \textit{``historic behavior"} in literature. Namely, \textit{an orbit keeps having new ideas about what it wants to do and keeps changing considerably so that the values of the time averages give information about the epoch to which the sojourn time belongs: it has a history}. Hence, \textit{the orbit witnesses the history of the dynamical system and records the fluctuations}.

In this paper, we propose two vast classes of replicator equations generated by \textit{similar-order preserving mappings} which exhibit \textit{stable} as well as \textit{mean historic behavior}. In the latter case, the time averages of the orbit will slowly oscillate during the evolution of the system and do not converge to any limit. This will eventually cause the divergence of higher-order repeated time averages.

\section{The Main Results}

Let $\mathbf{I}_m:=\{1,\cdots, m\}$ be a finite set and $\{\mathbf{e}_k\}_{k=1}^m$ be
the standard basis of $\mathbb{R}^m$. Suppose that $\mathbb{R}^m$ is equipped with the
$l_1-$norm  $\|\mathbf{x}\|_1:=\sum_{k=1}^{m}|x_k|$ where $\mathbf{x}:=(x_1,\cdots,x_m) \in \mathbb{R}^m$. We also set $\langle\mathbf{x},\mathbf{y}\rangle:=\sum_{k=1}^{m}x_ky_k$ for any two vectors $\mathbf{x},\mathbf{y}\in\mathbb{R}^{m}$.
We say that $\mathbf{x}\geq 0$  if $x_{k}\geq 0$ for all $k\in\mathbf{I}_m$. Moreover, we say that $\mathbf{x}\geq \mathbf{y}$  if one has $\mathbf{x}-\mathbf{y} \geq 0$, i.e., $x_{k}\geq y_{k}$ for all $k\in\mathbf{I}_m$. Let $\mathbb{R}_{+}^m:=\{\mathbf{x}\in\mathbb{R}^{m}: \mathbf{x}\geq 0\}$ be a positive orthant, $\mathbb{B}_{+}^m:=\{\mathbf{x}\in\mathbb{R}_{+}^{m}: \|\mathbf{x}\|_1\leq 1\}$ be a part of unit ball located in the positive orthant, and $\mathbb{S}^{m-1}:=\{\mathbf{x}\in\mathbb{R}_{+}^{m}: \|\mathbf{x}\|_1=1\}$ be the standard simplex. For a vector 
$\mathbf{x}\in\mathbb{R}_{+}^m$, we set 
$$\textup{\textbf{supp}}(\mathbf{x})=\{k\in \mathbf{I}_m: x_k>0\} \quad \textup{and} \quad \textup{\textbf{null}}(\mathbf{x})=\{k\in\mathbf{I}_m: x_k=0\}.$$  

Two vectors $\mathbf{x},\mathbf{y}\in\mathbb{R}_{+}^{m}$ are called \textit{similarly ordered}, denoted by $\mathbf{x}\approx\mathbf{y}$, if one has that $x_i\left(\gtreqqless\right)x_j$ if and only if $y_i\left(\gtreqqless\right)y_j$ for all $i,j\in\mathbf{I}_m$. Obviously, for a vector $\mathbf{x}\in\mathbb{R}_{+}^{m}$,  a set $\textsc{\textbf{SOC}}[\mathbf{x}]=\{\mathbf{y}\in\mathbb{R}_{+}^{m}: \mathbf{x}\approx\mathbf{y}\}$ of all vectors $\mathbf{y}\in\mathbb{R}_{+}^{m}$ for which $\mathbf{x}$ and $\mathbf{y}$ are similarly ordered  is a \textit{convex cone}. 
\begin{definition}
	A continuous mapping $\mathbf{F}:\mathbb{R}^m_{+}\to\mathbb{R}^m_{+}$ is called \textit{similar-order preserving} if one has $\mathbf{F}(\mathbf{x})\in\textsc{\textbf{SOC}}[\mathbf{x}]$, i.e, $\mathbf{F}(\mathbf{x})\approx\mathbf{x}$ for all $\mathbf{x}\in\mathbb{R}_{+}^{m}$.
\end{definition}

In the next section, we provide some examples for similar-order preserving mappings. 

Let $\mathbf{F}:\mathbb{B}^m_{+}\to\mathbb{R}^m_{+}$, $\mathbf{F}(\mathbf{x}):=\left(f_1(\mathbf{x}),\cdots,f_m(\mathbf{x})\right)$ be a continuously differentiable similar-order preserving mapping. Throughout this paper, without loss of generality,  we may assume that $0<\mathbf{F}(\mathbf{x})\leq \mathbbm{1}$ for all $\mathbf{x}\in\mathbb{B}_{+}^{m}$ where $\mathbbm{1}=(1,1,\cdots,1)$. Otherwise, for small $\varepsilon>0$, we consider a continuous mapping $\widehat{\mathbf{F}}:\mathbb{B}^m_{+}\to\mathbb{R}^m_{+}$,  $$\widehat{\mathbf{F}}(\mathbf{x}):=\frac{1}{M_1}\left(\mathbf{F}(\mathbf{x})+\varepsilon \mathbbm{1}\right)$$ which is similar-order preserving and $0<\widehat{\mathbf{F}}(\mathbf{x})\leq \mathbbm{1}$ for all $\mathbf{x}\in\mathbb{B}_{+}^{m}$ where  
$$M_1:=\max\{M+\varepsilon,1\} \quad \textup{and} \quad M:=\max\limits_{k\in\mathbf{I}_m}\max\limits_{\mathbf{x}\in\mathbb{B}_{+}^m}f_k(\mathbf{x}).$$

We consider a discrete-time replicator equation  $\mathcal{R}_{S}:\mathbb{S}^{m-1}\to\mathbb{S}^{m-1}$ generated by a continuously differentiable similar-order preserving mapping $\mathbf{F}:\mathbb{B}^m_{+}\to\mathbb{R}^m_{+}$,  $\mathbf{F}(\mathbf{x}):=\left(f_1(\mathbf{x}),\cdots,f_m(\mathbf{x})\right)$ 
\begin{equation}\label{RES}
	\left(\mathcal{R}_{S}(\mathbf{x})\right)_k=x_k\left(1+f_k(\mathbf{x})-\sum\limits_{i=1}^mx_if_i(\mathbf{x})\right), \quad \forall \ k\in\mathbf{I}_m.
\end{equation}

Obviously, in order the simplex $\mathbb{S}^{m-1}$ to be an invariant set for the replicator equation given by \eqref{RES}, we must have some constraints on a continuous similar-order preserving mapping $\mathbf{F}:\mathbb{B}^m_{+}\to\mathbb{R}^m_{+}$,  $\mathbf{F}(\mathbf{x}):=\left(f_1(\mathbf{x}),\cdots,f_m(\mathbf{x})\right)$. For example, the following constraint
$$0\leq \max_{k\in\mathbf{I}_m}\max_{\mathbf{x}\in\mathbb{S}^{m-1}}f_k(\mathbf{x})-\min_{k\in\mathbf{I}_m}\min_{\mathbf{x}\in\mathbb{S}^{m-1}}f_k(\mathbf{x})\leq 1$$  is sufficient for the simplex $\mathbb{S}^{m-1}$ to be invariant. Indeed, since $$\min_{k\in\mathbf{I}_m}\min_{\mathbf{x}\in\mathbb{S}^{m-1}}f_k(\mathbf{x})\leq\sum\limits_{i=1}^mx_if_i(\mathbf{x})\leq\max_{k\in\mathbf{I}_m}\max_{\mathbf{x}\in\mathbb{S}^{m-1}}f_k(\mathbf{x}),\quad \forall \ \mathbf{x}\in\mathbb{S}^{m-1},$$ we obtain for all $\mathbf{x}\in\mathbb{S}^{m-1}$ and $k\in\mathbf{I}_m$ that 
\begin{multline*}
	1+f_k(\mathbf{x})-\sum\limits_{i=1}^mx_if_i(\mathbf{x})\geq 1+f_k(\mathbf{x})-\max_{k\in\mathbf{I}_m}\max_{\mathbf{x}\in\mathbb{S}^{m-1}}f_k(\mathbf{x})\geq\\ f_k(\mathbf{x})-\min_{k\in\mathbf{I}_m}\min_{\mathbf{x}\in\mathbb{S}^{m-1}}f_k(\mathbf{x})\geq 0.
\end{multline*}
Hence, we have $\sum_{k=1}^m\left(\mathcal{R}_{S}(\mathbf{x})\right)_k=\sum_{k=1}^mx_k=1$ and $\left(\mathcal{R}_{S}(\mathbf{x})\right)_k\geq 0$ for all $\mathbf{x}\in\mathbb{S}^{m-1}$ and $k\in\mathbf{I}_m$. Since, throughout this paper, we always assume that $0<\mathbf{F}(\mathbf{x})\leq \mathbbm{1}$ for all $\mathbf{x}\in\mathbb{B}_{+}^{m}$ where $\mathbbm{1}=(1,1,\cdots,1)$, it is needless to say that the constraint considered above
$$0\leq \max_{k\in\mathbf{I}_m}\max_{\mathbf{x}\in\mathbb{S}^{m-1}}f_k(\mathbf{x})-\min_{k\in\mathbf{I}_m}\min_{\mathbf{x}\in\mathbb{S}^{m-1}}f_k(\mathbf{x})\leq 1$$ is always satisfied.

Recall (see \cite{Cressman2003,HofbauerSigmund2003,Sandholm,Sigmund2010}) that $\mathbf{x}\in\mathbb{S}^{m-1}$ is called a \textit{Nash equilibrium} of the discrete-time replicator equation  given by \eqref{RES} if $\langle\mathbf{x},\mathbf{F}(\mathbf{x})\rangle\geq \langle\mathbf{y},\mathbf{F}(\mathbf{x})\rangle$ for all $\mathbf{y}\in\mathbb{S}^{m-1}$. Moreover, $\mathbf{x}\in\mathbb{S}^{m-1}$ is called a \textit{strictly Nash equilibrium} if $\langle\mathbf{x},\mathbf{F}(\mathbf{x})\rangle > \langle\mathbf{y},\mathbf{F}(\mathbf{x})\rangle$ for all $\mathbf{y}\in\mathbb{S}^{m-1}$ but $\mathbf{y}\neq\mathbf{x}$. A point $\mathbf{x}\in\mathbb{S}^{m-1}$ is called \textit{a rest (fixed) point} if $\mathcal{R}_{S}(\mathbf{x})=\mathbf{x}$. A rest (fixed) point $\mathbf{x}\in\mathbb{S}^{m-1}$ is called \textit{stable} if for every neighborhood $U(\mathbf{x})\subset\mathbb{S}^{m-1}$ of $\mathbf{x}$ there exists another neighborhood $V(\mathbf{x})\subset U(\mathbf{x})\subset\mathbb{S}^{m-1}$ of $\mathbf{x}$ such that an orbit $\{\mathbf{y}, \mathcal{R}_{S}(\mathbf{y}), \cdots, \mathcal{R}_{S}^{(n)}(\mathbf{y}), \cdots\}$ of any initial point $\mathbf{y}\in V(\mathbf{x})$ remains inside of the neighborhood $U(\mathbf{x})$. A rest (fixed) point $\mathbf{x}\in\mathbb{S}^{m-1}$ is called \textit{attracting} if there exists a neighborhood $V(\mathbf{x})\subset\mathbb{S}^{m-1}$ of $\mathbf{x}$ such that an orbit $\{\mathbf{y}, \mathcal{R}_{S}(\mathbf{y}), \cdots, \mathcal{R}_{S}^{(n)}(\mathbf{y}), \cdots\}$ of any initial point $\mathbf{y}\in V(\mathbf{x})$ converges to $\mathbf{x}$. A rest (fixed) point $\mathbf{x}\in\mathbb{S}^{m-1}$ is called \textit{asymptotically stable} if it is both \textit{stable} and \textit{attracting}. 

The following result shows that a vast class of discrete-time replicator equations satisfies \textit{``The Folk Theorem of Evolutionary Game Theory"}.

\begin{atheorem}\label{FTEGT}
	The following statements are true  for the discrete-time replicator equation $\mathcal{R}_{S}:\mathbb{S}^{m-1}\to\mathbb{S}^{m-1}$ given by \eqref{RES}\textup{:}
	\begin{itemize}
		\item[$(i)$]  A Nash equilibrium is a rest point; 
		\item[$(ii)$] A stable rest point is a Nash equilibrium; 
		\item[$(iii)$] A strictly Nash equilibrium is asymptotically stable; 
		\item[$(iv)$] Any interior convergent orbit evolves to a Nash equilibrium.
	\end{itemize}
\end{atheorem}

It is worth mentioning that some particular cases of the discrete-time replicator equations $\mathcal{R}_{S}:\mathbb{S}^{m-1}\to\mathbb{S}^{m-1}$ given by \eqref{RES} were studied in the literature (see \cite{Saburov2021a,Saburov2022}).

We now present a discrete-time replicator equation of a zero-sum game that exhibits \textit{mean historic behavior}. We consider either one of the following discrete-time replicator equations $\mathcal{R}_{H}:\mathbb{S}^{2}\to\mathbb{S}^{2}$ of a zero-sum game defined by a continuously differentiable similar-order preserving mapping $\mathbf{F}:\mathbb{B}^3_{+}\to\mathbb{R}^3_{+}$,  $\mathbf{F}(\mathbf{x}):=\left(f_1(\mathbf{x}),f_2(\mathbf{x}),f_3(\mathbf{x})\right)$
{\begin{eqnarray}\label{REH1}
				\mathcal{R}_{H}:\left\{
				\begin{matrix}
						(\mathcal{R}_{H}(\mathbf{x}))_1=x_1\left(1+x_2f_1(\mathbf{x})-x_3f_3(\mathbf{x})\right)\\
						(\mathcal{R}_{H}(\mathbf{x}))_2=x_2\left(1+x_3f_2(\mathbf{x})-x_1f_1(\mathbf{x})\right)\\
						(\mathcal{R}_{H}(\mathbf{x}))_3=x_3\left(1+x_1f_3(\mathbf{x})-x_2f_2(\mathbf{x})\right)
					\end{matrix}
				\right.
			\end{eqnarray}
		or
		\begin{eqnarray}\label{REH2}
				\mathcal{R}_{H}:\left\{
				\begin{matrix}
						(\mathcal{R}_{H}(\mathbf{x}))_1=x_1\left(1+x_3f_1(\mathbf{x})-x_2f_2(\mathbf{x})\right)\\
						(\mathcal{R}_{H}(\mathbf{x}))_2=x_2\left(1+x_1f_2(\mathbf{x})-x_3f_3(\mathbf{x})\right)\\
						(\mathcal{R}_{H}(\mathbf{x}))_3=x_3\left(1+x_2f_3(\mathbf{x})-x_1f_1(\mathbf{x})\right)
					\end{matrix}
				\right.
		\end{eqnarray}}
		
		It is wroth mentioning that a \textit{density-dependent} payoff matrix of a zero-sum game is given by
		\begin{eqnarray*}
			\left(
			\begin{matrix}
				0 & f_1(\mathbf{x}) & -f_3(\mathbf{x})\\
				-f_1(\mathbf{x}) & 0 & f_2(\mathbf{x}) \\
				f_3(\mathbf{x}) & -f_2(\mathbf{x}) & 0
			\end{matrix}
			\right) 
			\quad \textup{or} \quad
			\left(
			\begin{matrix}
				0 & -f_2(\mathbf{x}) & f_1(\mathbf{x})\\
				f_2(\mathbf{x}) & 0 & -f_3(\mathbf{x}) \\
				-f_1(\mathbf{x}) & f_3(\mathbf{x}) & 0
			\end{matrix}
			\right).
		\end{eqnarray*}

		In the language of evolutionary game theory, the discrete-time replicator equation $\mathcal{R}_{H}:\mathbb{S}^{2}\to\mathbb{S}^{2}$ given by \eqref{REH1} or \eqref{REH2} represents the interactions of three competing species where \textit{the species 1 beats (loses to) the species 2, the species 2 beats (loses to) the species 3, and the species 3 beats (loses to) the species 1} (like the \textit{Rock-Paper-Scissors} game). Some particular cases of the discrete-time replicator equations $\mathcal{R}_{H}:\mathbb{S}^{2}\to\mathbb{S}^{2}$ were studied in the literature (see \cite{KBMM2010,GanikhodzhaevZanin,JamilovMukhamedov20222,JamilovMukhamedov20221,JamilovScheutzowVorkastner2019,Saburov2015a,Saburov2015b,Saburov2018c,Saburov2019a,Saburov2020,Saburov2021c,Saburov2021d,Vallander1972,Vallander2007,Zakharevich1978}). 
		
		We say that the discrete-time replicator equation $\mathcal{R}_{H}:\mathbb{S}^{2}\to\mathbb{S}^{2}$ given by \eqref{REH1} or \eqref{REH2} has \textit{mean historic behavior} if the set of initial points $\mathbf{x}\in\mathbb{S}^{2}$ which give rise to orbits with divergent time averages $\frac{1}{n}\sum_{k=0}^{n-1}\mathcal{R}_{H}^{(k)}(\mathbf{x})$ has positive Lebesgue measure (see \cite{Saburov2021d}). We define the $s^{th}-$order repeated time averages $\left\{\mathcal{A}^{(s)}_{n}\left(\mathbf{x}\right)\right\}_{n=1}^\infty$ for all $s\in\mathbb{N}$ as follows
		$$\mathcal{A}^{(s)}_{n}\left(\mathbf{x}\right):=\frac{1}{n}\sum_{k=1}^{n}\mathcal{A}^{(s-1)}_{k}\left(\mathbf{x}\right), \quad \mathcal{A}^{(1)}_{n}\left(\mathbf{x}\right):=\frac{1}{n}\sum_{k=0}^{n-1}\mathcal{R}_{H}^{(k)}(\mathbf{x}).$$
		
		\begin{atheorem}\label{HistircBehavior}
			The discrete-time replicator equation $\mathcal{R}_{H}:\mathbb{S}^{2}\to\mathbb{S}^{2}$ given by \eqref{REH1} or \eqref{REH2} has mean historic behavior. Moreover, for any interior initial point, the $s^{th}-$order repeated time averages $\left\{\mathcal{A}^{(s)}_{n}\left(\mathbf{x}\right)\right\}_{n=1}^\infty$ do not converge for any $s\in\mathbb{N}$.
		\end{atheorem}

		\section{Some Examples for Similar-Order Preserving Mappings}
		
		Let $\mathbf{F},\mathbf{G}:\mathbb{R}^m_{+}\to\mathbb{R}^m_{+}$, $\mathbf{F}(\mathbf{x}):=\left(f_1(\mathbf{x}),\cdots,f_m(\mathbf{x})\right)$,  $\mathbf{G}(\mathbf{x}):=\left(g_1(\mathbf{x}),\cdots,g_m(\mathbf{x})\right)$  be continuous similar-order preserving mappings and $\varphi,\psi:\mathbb{R}^m_{+}\to (0,\infty)$  be continuous functions.  Let $h:\mathbb{R}_{+}\to\mathbb{R}_{+}$ be a strictly increasing continuous function and $\mathbf{H}:\mathbb{R}^m_{+}\to\mathbb{R}^m_{+}$, $\mathbf{H}(\mathbf{x}):=\left(h(x_1),\cdots,h(x_m)\right)$. It is easy to check that the following  mappings
	
	{\small	\begin{itemize}
			\item[$(i)$] $\mathbf{F}_{\varphi,\psi}:\mathbb{R}^m_{+}\to\mathbb{R}^m_{+}$ where $\mathbf{F}_{\varphi,\psi}(\mathbf{x}):=\bigl(\varphi(\mathbf{x})f_1(\mathbf{x})+\psi(\mathbf{x}),\cdots,\varphi(\mathbf{x})f_m(\mathbf{x})+\psi(\mathbf{x})\bigr)$;
			\item[$(ii)$] $\varphi\mathbf{F}+\psi\mathbf{G}:\mathbb{R}^m_{+}\to\mathbb{R}^m_{+}$ where  $\left(\varphi\mathbf{F}+\psi\mathbf{G}\right)(\mathbf{x}):=\varphi(\mathbf{x})\mathbf{F}(\mathbf{x})+\psi(\mathbf{x})\mathbf{G}(\mathbf{x})$;
			\item[$(iii)$] $\mathbf{H}\circ\mathbf{F}:\mathbb{R}^m_{+}\to\mathbb{R}^m_{+}$ where $\left(\mathbf{H}\circ\mathbf{F}\right)(\mathbf{x}):=\mathbf{H}\left(\mathbf{F}(\mathbf{x})\right)=\bigl(h(f_1(\mathbf{x})),\cdots,h(f_m(\mathbf{x}))\bigr)$;
			\item[$(iv)$] $\mathbf{F}\circ\mathbf{H}:\mathbb{R}^m_{+}\to\mathbb{R}^m_{+}$ where $\left(\mathbf{F}\circ\mathbf{H}\right)(\mathbf{x}):=\mathbf{F}\left(\mathbf{H}(\mathbf{x})\right)=\bigl(f_1(\mathbf{H}(\mathbf{x})),\cdots,f_m(\mathbf{H}(\mathbf{x}))\bigr)$;
			\item[$(v)$] $\mathbf{F}\circ\mathbf{G}:\mathbb{R}^m_{+}\to\mathbb{R}^m_{+}$ where $\left(\mathbf{F}\circ\mathbf{G}\right)(\mathbf{x}):=\mathbf{F}\left(\mathbf{G}(\mathbf{x})\right)=\bigl(f_1(\mathbf{G}(\mathbf{x})),\cdots,f_m(\mathbf{G}(\mathbf{x}))\bigr)$ 
		\end{itemize} } \noindent are also similar-order preserving. Moreover, if $\mathbf{F}(\mathbf{x})>0$ and $\mathbf{G}(\mathbf{x})>0$ for all $\mathbf{x}\in\mathbb{R}_{+}^{m}$ then the mapping
		{\small \begin{itemize}
			\item[$(vi)$] $\mathbf{F}\odot\mathbf{G}:\mathbb{R}^m_{+}\to\mathbb{R}^m_{+}$ where $\left(\mathbf{F}\odot\mathbf{G}\right)(\mathbf{x}):=\bigl(f_1(\mathbf{x})\cdot g_1(\mathbf{x}),\cdots,f_m(\mathbf{x})\cdot g_m(\mathbf{x})\bigr)$ 
		\end{itemize} }
		\noindent is also similar-order preserving. 
		
		These properties allow to generate various kinds of similar-order preserving mappings. Let us now provide some concrete examples for similar-order preserving mappings. 
		
		Let $x_{[1]}\geq \cdots\geq x_{[m]}$ denote the components of $\textbf{x}=(x_1,\cdots,x_m)\in\mathbb{R}^m$ in a non-increasing order. We say that $\mathbf{x}$ \textit{majorizes} $\mathbf{y}$ (or $\mathbf{y}$ is \textit{majorized} by $\mathbf{x}$), written $\mathbf{x}\succ \mathbf{y}$, if one has $
		\sum_{i=1}^{k}x_{[i]}\geq \sum_{i=1}^{k}y_{[i]}$ for all $k\in\mathbf{I}_{m-1}$ and $\sum_{i=1}^{m}x_{i}=\sum_{i=1}^{m}y_{i}$. 
		
		Recall (see \cite{MOA}) that a function $\varphi:\mathbb{B}_{+}^m\to\mathbb{R}$ is said to be (\textit{strictly}) \textit{increasing} if $\mathbf{x}\geq\mathbf{y}$ (resp. $\mathbf{x}\geq\mathbf{y}$ but $\mathbf{x}\neq\mathbf{y}$) implies $\varphi(\mathbf{x})\geq\varphi(\mathbf{y})$ (resp. $\varphi(\mathbf{x})>\varphi(\mathbf{y})$). A function $\varphi:\mathbb{B}_{+}^m\to\mathbb{R}$ is said to be (\textit{strictly}) \textit{Schur-convex} if  $\mathbf{x}\succ\mathbf{y}$ (resp. $\mathbf{x}\succ\mathbf{y}$ but $\mathbf{y}$ is not a permutation
		of $\mathbf{x}$) implies $\varphi(\mathbf{x})\geq\varphi(\mathbf{y})$ (resp. $\varphi(\mathbf{x})>\varphi(\mathbf{y})$).
		
		It is worth mentioning (see \cite{MOA}) that, since $\varphi:\mathbb{B}_{+}^m\to\mathbb{R}$ is a strictly Schur-convex function, we have that $\frac{\partial\varphi}{\partial x_i}(\mathbf{x}) \left(\gtrless\right) \frac{\partial\varphi}{\partial x_j}(\mathbf{x})$ if and only if $x_i\left(\gtrless\right)x_j$ for $i,j\in\mathbf{I}_m$. Moreover, since $\varphi:\mathbb{B}_{+}^m\to\mathbb{R}$ is symmetric,  we also have that $\frac{\partial\varphi}{\partial x_i}(\mathbf{x}) = \frac{\partial\varphi}{\partial x_j}(\mathbf{x})$ if and only if $x_i=x_j$ for $i,j\in\mathbf{I}_m$. 
		
		Hence, for any symmetric, continuously differentiable, strictly increasing and strictly Schur-convex function $\varphi:\mathbb{R}_{+}^m\to\mathbb{R}_{+}$, a gradient vector field $\nabla\varphi:\mathbb{R}_{+}^m\to\mathbb{R}_{+}^m$ where $$\nabla \varphi(\mathbf{x}):=\left(\frac{\partial\varphi}{\partial x_1}(\mathbf{x}),\cdots,\frac{\partial\varphi}{\partial x_m}(\mathbf{x})\right)$$ is always similar-order preserving. Particularly, the gradient vector fields of the following functions are also similar-order preserving (the reader may refer to the monograph \cite{MOA}):
		\begin{itemize}
			\item[$(i)$] \textit{Complete Symmetric Functions}: $\varphi_k:\mathbb{R}_{+}^m\to\mathbb{R}_{+}$, $$\varphi_{k}(\mathbf{x})=\sum\limits_{i_1+i_2+\cdots+i_m=k}x_1^{i_1}x_2^{i_2}\cdots x_m^{i_m}, \qquad \forall \  k\in\mathbf{I}_m;$$ 
			\item[$(ii)$] \textit{Symmetric Gauge Functions}: $\varphi_p:\mathbb{R}_{+}^m\to\mathbb{R}_{+}$,
			$$\varphi_p(\mathbf{x})=\|\mathbf{x}\|_p=\sqrt[p]{\sum_{k=1}^{m}x^{p}_k}, \qquad p>1;$$
			\item[$(iii)$] \textit{Gamma Functions}: $\varphi:\mathbb{R}_{+}^m\to\mathbb{R}_{+}$,
			$$\varphi(\mathbf{x})=\prod_{k=1}^{m}\Gamma(x_k+a), \qquad a\geq 1,$$ where $\Gamma:\mathbb{R}_{+}\to\mathbb{R}_{+}$ is a gamma function;
			\item[$(iv)$] \textit{Convex Separate Variables Functions}: $\varphi:\mathbb{R}_{+}^m\to\mathbb{R}_{+}$,
			$$\varphi(\mathbf{x})=\sum_{k=1}^{m}f(x_k),$$ where $f:\mathbb{R}_{+}\to\mathbb{R}_{+}$ is a continuously differentiable, strictly increasing, and strictly convex function; 
			\item[$(v)$] \textit{Symmetric Composite Functions}: $\varphi_k:\mathbb{R}_{+}^m\to\mathbb{R}_{+}$,
			$$\varphi_k(\mathbf{x})=\sum\limits_{1\leq i_1<\cdots<i_k\leq m}f(x_{i_1})\cdots f(x_{i_k}), \qquad \forall \  k\in\mathbf{I}_m,$$ where $f:\mathbb{R}_{+}\to\mathbb{R}_{+}$ is a continuously differentiable, strictly increasing, and strictly \textit{log-convex} function;
			\item[$(vi)$] \textit{General Convex Symmetric Functions}: $\varphi:\mathbb{R}_{+}^m\to\mathbb{R}_{+}$ is a continuously differentiable, symmetric, strictly increasing, and strictly convex function.
			\item[$(vii)$] \textit{Composition Functions}: If $\psi:\mathbb{R}_{+}^m\to\mathbb{R}_{+}$ is a continuously differentiable, strictly increasing, and strictly Schur-convex function and $h:\mathbb{R}_{+}\to\mathbb{R}_{+}$ is a continuously differentiable, strictly increasing and strictly convex function then $\varphi:\mathbb{R}_{+}^m\to\mathbb{R}_{+}$, $\varphi(\mathbf{x})=\psi(h(x_1),\cdots, h(x_m))$ is also a continuously differentiable, strictly increasing, and strictly Schur-convex function.
		\end{itemize}

		\section{The Proof of Theorem \ref{FTEGT}} 
		
		In order to study the stability of rest (fixed) points of the discrete-time replicator equations $\mathcal{R}_{S}:\mathbb{S}^{m-1}\to\mathbb{S}^{m-1}$ given by \eqref{RES}, we employ a Lyapunov function. Recall (see \cite{Cressman2003,HofbauerSigmund1998,Sandholm,Sigmund2010}) that a continuous function $\varphi:\mathbb{S}^{m-1}\to\mathbb{R}$ is called a \textit{Lyapunov function} if the number sequence $\{\varphi(\mathbf{x}), \varphi(\mathcal{R}_{S}(\mathbf{x})), \cdots, \varphi(\mathcal{R}_{S}^{(n)}(\mathbf{x})), \cdots \}$ is a bounded monotone sequence for any initial point $\mathbf{x}\in\mathbb{S}^{m-1}$.

		\begin{description}
			\item[Part $(i)$] Let us first show that the center $\mathbf{c}_\alpha:=\frac{1}{|\alpha|}\sum_{k\in\alpha}\mathbf{e}_k$ of any face $$\mathbb{S}_{\alpha}:=\{\mathbf{x}\in\mathbb{S}^{m-1}:\textup{\textbf{supp}}(\mathbf{x})\subset \alpha\}$$ of the simplex $\mathbb{S}^{m-1}$ is a rest (fixed) point for all $\alpha\subset \mathbf{I}_m$. Indeed, it is clear that if $k\not\in\alpha$ then  $(\mathcal{R}_S(\mathbf{c}_\alpha))_k=(\mathbf{c}_\alpha)_k=0$ and if  $k\in\alpha$ then, since $\mathbf{F}:\mathbb{B}^m_{+}\to\mathbb{R}^m_{+}$ is a similar-order preserving mapping, we have $f_k(\mathbf{c}_\alpha)=\sum\limits_{i\in\alpha}(\mathbf{c}_\alpha)_if_i(\mathbf{c}_\alpha)$ and $(\mathcal{R}_S(\mathbf{c}_\alpha))_k=(\mathbf{c}_\alpha)_k=\frac{1}{|\alpha|}$.
			
			Moreover, there is no any other rest (fixed) point except the center of every face of the simplex. Indeed, let $\mathbf{x}\in\mathbb{S}^{m-1}$ be any rest (fixed) point. We then obtain  $(\mathcal{R}_S(\mathbf{x}))_k=x_k$ for all $k\in\textup{\textbf{supp}}(\mathbf{x})$. This is equivalent saying that $f_{k_1}(\mathbf{x})=\sum\limits_{i\in\alpha}x_if_i(\mathbf{x})=f_{k_2}(\mathbf{x})$ for any $k_1,k_2\in\textup{\textbf{supp}}(\mathbf{x})$ and $k_1\neq k_2$. Since $\mathbf{F}:\mathbb{B}^m_{+}\to\mathbb{R}^m_{+}$ is a similar-order preserving mapping, we obtain $x_{k_1}=x_{k_2}$ for any $k_1,k_2\in\textup{\textbf{supp}}(\mathbf{x})$ and $k_1\neq k_2$. This means that $\mathbf{x}$ is the center of the face $\mathbb{S}_{\alpha}$ where $\alpha=\textup{\textbf{supp}}(\mathbf{x})$.

			We now describe all Nash equilibria. Since
			$$\langle\mathbf{c}_\alpha,\mathbf{F}(\mathbf{c}_\alpha)\rangle=\max\limits_{k\in\mathbf{I}_m}f_k(\mathbf{c}_\alpha) \geq \langle\mathbf{y},\mathbf{F}(\mathbf{c}_\alpha)\rangle$$ for all $\mathbf{y}\in\mathbb{S}^{m-1},$ the center of every face of the simplex is a Nash equilibrium. 
			
			Moreover, if $\mathbf{x}\in\mathbb{S}^{m-1}$ is a Nash equilibrium then $$\max\limits_{k\in\mathbf{I}_m}f_k(\mathbf{x}) \geq\langle\mathbf{x},\mathbf{F}(\mathbf{x})\rangle\geq \langle\mathbf{e}_j,\mathbf{F}(\mathbf{x})\rangle=\max\limits_{k\in\mathbf{I}_m}f_k(\mathbf{x})$$ 
			where $j\in\mathbf{I}_m$ is an index such that $f_j(\mathbf{x})=\max\limits_{k\in\mathbf{I}_m}f_k(\mathbf{x})$. Since 
			$$f_i(\mathbf{x}) \left(\gtreqqless\right) f_j(\mathbf{x})\quad \textup{if and only if} \quad x_i\left(\gtreqqless\right)x_j \quad \textup{for all} \quad i,j\in\mathbf{I}_m,$$ we then obtain $\mathbf{x}=\mathbf{c}_\alpha$ for a subset $$\alpha:=\left\{j\in\mathbf{I}_m: f_j(\mathbf{x})=\max\limits_{k\in\mathbf{I}_m}f_k(\mathbf{x})\right\}\subset\mathbf{I}_m.$$ This shows that every Nash equilibrium is a rest (fixed) point.
			
			\item[Part $(ii)$] Since every face $\mathbb{S}_{\alpha}$ is invariant, we first show that $$\mathcal{M}_{\alpha, k}(\mathbf{x}):=\max\limits_{i\in\alpha}x_{i}-x_k$$ is an increasing Lyapunov function over the interior $$\textup{int}\mathbb{S}_{\alpha}:=\{\mathbf{x}\in\mathbb{S}_{\alpha}:\textup{\textbf{supp}}(\mathbf{x})=\alpha\}$$ of the face $\mathbb{S}_{\alpha}$ for every $\alpha\subset\mathbf{I}_m$ and for every $k\in\alpha$. Indeed, for any $k,t\in\alpha$, if $x_t=x_k$ then $(\mathcal{R}_S(\mathbf{x}))_t=(\mathcal{R}_S(\mathbf{x}))_k$ and if $x_t\neq x_k$ then 
			$$\frac{(\mathcal{R}_S(\mathbf{x}))_t-(\mathcal{R}_S(\mathbf{x}))_k}{x_t-x_k}=\left[\frac{x_tf_t(\mathbf{x})-x_kf_k(\mathbf{x})}{x_t-x_k}+\sum\limits_{i\in\alpha}x_i\left(1-f_i(\mathbf{x})\right)\right].
			$$
			
			Since $0 < \mathbf{F}(\mathbf{x})\leq \mathbbm{1}$ and $x_tf_t(\mathbf{x}) \left(\gtreqqless\right) x_kf_k(\mathbf{x})$ if and only if $x_t\left(\gtreqqless\right)x_k$ for all $k,t\in\alpha$, we obtain that $(\mathcal{R}_S(\mathbf{x}))_t \left(\gtreqqless\right) (\mathcal{R}_S(\mathbf{x}))_k$ if and only if $x_t\left(\gtreqqless\right)x_k$ for all $k,t\in\alpha$. This yields $$\mathsf{\textbf{MaxInd}}_\alpha\left(\mathcal{R}_S(\mathbf{x})\right)=\mathsf{\textbf{MaxInd}}_\alpha(\mathbf{x})$$ for any $\mathbf{x}\in\textup{int}\mathbb{S}_{\alpha}$ where $\mathsf{\textbf{MaxInd}}_\alpha(\mathbf{x}):=\{k\in\alpha: x_k=\max\limits_{i\in\alpha}x_{i}\}$. 
			
			Moreover, if $t \in \mathsf{\textbf{MaxInd}}_\alpha\left(\mathcal{R}_S(\mathbf{x})\right)=\mathsf{\textbf{MaxInd}}_\alpha(\mathbf{x})$ then for all $k\in\alpha$ with $x_t\neq x_k$ we have that 
			$$\frac{(\mathcal{R}_S(\mathbf{x}))_t-(\mathcal{R}_S(\mathbf{x}))_k}{x_t-x_k}=1+\left[\frac{x_k\left(f_t(\mathbf{x})-f_k(\mathbf{x})\right)}{x_t-x_k}+\sum\limits_{i\in\alpha}x_i\left(f_t(\mathbf{x})-f_i(\mathbf{x})\right)\right].$$
			
			This yields $\mathcal{M}_{\alpha, k}(\mathcal{R}_S(\mathbf{x}))\geq\mathcal{M}_{\alpha, k}(\mathbf{x})$ for all $\mathbf{x}\in\textup{int}\mathbb{S}_{\alpha}$ and for all $k\in\alpha$. By repeating this process, we obtain  $$\mathcal{M}_{\alpha, k}(\mathcal{R}^{(n+1)}_S(\mathbf{x}))\geq\mathcal{M}_{\alpha, k}(\mathcal{R}^{(n)}_S(\mathbf{x})), \ \mathsf{\textbf{MaxInd}}_\alpha\left(\mathcal{R}^{(n)}_S(\mathbf{x})\right)=\mathsf{\textbf{MaxInd}}_\alpha(\mathbf{x})$$ for all $\mathbf{x}\in\textup{int}\mathbb{S}_{\alpha}$, $k\in\alpha$ and $n\in\mathbb{N}$. 
			This shows that $\mathcal{M}_{\alpha, k}(\mathbf{x})$ is an increasing Lyapunov function. 
			
			We now show that an orbit starting from any initial point $\mathbf{x}\in\textup{int}\mathbb{S}_{\alpha}$ converges to the center $\mathbf{c}_\beta$ of the face $\mathbb{S}_{\beta}\subset\mathbb{S}_{\alpha}$ where $\beta=\mathsf{\textbf{MaxInd}}_\alpha(\mathbf{x})$. Indeed, for every $k\in\alpha$, we have 
			\begin{multline*}
				(\mathcal{R}^{(n)}_S(\mathbf{x}))_k=\frac{1}{|\alpha|}\left[1+\sum_{i\in\alpha}\left((\mathcal{R}^{(n)}_S(\mathbf{x}))_k-(\mathcal{R}^{(n)}_S(\mathbf{x}))_i\right)\right]=\\
				\frac{1}{|\alpha|}\left[1+\sum_{i\in\alpha}\left(\mathcal{M}_{\alpha, i}(\mathcal{R}^{(n)}_S(\mathbf{x}))-\mathcal{M}_{\alpha, k}(\mathcal{R}^{(n)}_S(\mathbf{x}))\right)\right].
			\end{multline*}
			Since the sequence $\left\{\mathcal{M}_{\alpha, k}(\mathcal{R}^{(n)}_S(\mathbf{x}))\right\}_{n=0}^\infty$ is convergent for every $k\in\alpha$, so is the sequence $\left\{(\mathcal{R}^{(n)}_S(\mathbf{x}))_k\right\}_{n=0}^\infty$ for every $k\in\alpha$. Consequently, an orbit $\left\{\mathcal{R}^{(n)}_S(\mathbf{x})\right\}_{n=0}^\infty$ is also convergent and its limit is a rest point that is the center $\mathbf{c}_\beta$ of a face $\mathbb{S}_{\beta}\subset\mathbb{S}_{\alpha}$ for some  $\beta\subset\alpha$. Since $\mathsf{\textbf{MaxInd}}_\alpha\left(\mathcal{R}^{(n)}_S(\mathbf{x})\right)=\mathsf{\textbf{MaxInd}}_\alpha(\mathbf{x})$, we have $\beta=\mathsf{\textbf{MaxInd}}_\alpha(\mathbf{c}_\beta)\supset\mathsf{\textbf{MaxInd}}_\alpha(\mathbf{x})$. Since $\mathcal{M}_{\alpha, k}(\mathcal{R}^{(n+1)}_S(\mathbf{x}))\geq\mathcal{M}_{\alpha, k}(\mathcal{R}^{(n)}_S(\mathbf{x}))$, we have $\beta=\mathsf{\textbf{MaxInd}}_\alpha(\mathbf{c}_\beta)\subset\mathsf{\textbf{MaxInd}}_\alpha(\mathbf{x})$. Hence, we obtain  $\beta=\mathsf{\textbf{MaxInd}}_\alpha(\mathbf{x})$.
			
			We are now ready to show that the stable rest points are only vertices of the simplex which are Nash equilibria.
			
			On the one hand, the center $\mathbf{c}_\alpha$ of the face $\mathbb{S}_{\alpha}$ is not stable for all $\alpha\subset\mathbf{I}_m$ with $|\alpha|\geq2$. Indeed, for any small neighborhood $U(\mathbf{c}_\alpha)$ of the center $\mathbf{c}_\alpha$ there are points $\mathbf{x},\mathbf{y}\in U(\mathbf{c}_\alpha)$ such that $|\mathsf{\textbf{MaxInd}}_\alpha(\mathbf{x})|=|\mathsf{\textbf{MaxInd}}_\alpha(\mathbf{y})|=1$ and $\mathsf{\textbf{MaxInd}}_\alpha(\mathbf{x})\neq\mathsf{\textbf{MaxInd}}_\alpha(\mathbf{y})$ for which their orbits converge to two different vertices of the simplex.
			
			On the other hand, the vertex $\mathbf{e}_k$ of the simplex is stable for all $k\in\mathbf{I}_m$. Indeed, for any small neighborhood $U(\mathbf{e}_k)$ of the vertex $\mathbf{e}_k$ one has $$\mathsf{\textbf{MaxInd}}\left(\mathcal{R}_S(\mathbf{x})\right)=\mathsf{\textbf{MaxInd}}(\mathbf{x})=\{k\}$$ for all $\mathbf{x}\in U(\mathbf{e}_k)$ and
			$$
			\left(\mathcal{R}_{S}(\mathbf{x})\right)_k=x_k\left[1+\sum\limits_{i=1}^mx_i\left(f_k(\mathbf{x})-f_i(\mathbf{x})\right)\right]\geq x_k.
			$$
			Hence, we obtain  $\|\mathcal{R}_S(\mathbf{x})-\mathbf{e}_k\|_1=2(1-\left(\mathcal{R}_S(\mathbf{x})\right)_k)\leq 2(1-x_k)=\|\mathbf{x}-\mathbf{e}_k\|_1$ which implies $\mathcal{R}_S(U(\mathbf{e}_k))\subset U(\mathbf{e}_k)$ and consequently $\mathcal{R}_S^{(n)}(U(\mathbf{e}_k))\subset U(\mathbf{e}_k)$. This means that the stable rest points are only vertices of the simplex which are Nash equilibria.
			
			\item[Part $(iii)$] We first show that strictly Nash equilibria are only vertices of the simplex. 
			
			On the one hand, since $f_k(\mathbf{e}_k)>f_i(\mathbf{e}_k)$ for all $i\neq k$, we have that 
			$$\langle\mathbf{e}_k,\mathbf{F}(\mathbf{e}_k)\rangle=f_k(\mathbf{e}_k)=\max\limits_{i\in\mathbf{I}_m}f_i(\mathbf{e}_k) >\langle\mathbf{x},\mathbf{F}(\mathbf{e}_k)\rangle=\sum_{i\in\textup{\textbf{supp}}(\mathbf{x})}x_if_i(\mathbf{e}_k)$$ for all $\mathbf{x}\in {\mathbb{S}}^{m-1}\setminus\{\mathbf{e}_k\}$. 
			This shows that every vertex $\mathbf{e}_k$ of the simplex is a strictly Nash equilibrium for all $k\in\mathbf{I}_m$. 
			
			On the other hand, for any $\alpha\subset \mathbf{I}_m$ with $|\alpha|\geq2$ and $k\in\alpha$, we have that  
			$$ \langle\mathbf{c}_\alpha,\mathbf{F}(\mathbf{c}_\alpha)\rangle=\left(\sum_{i\in\alpha}\frac{1}{|\alpha|}\right)\max\limits_{i\in\alpha}f_i(\mathbf{c}_\alpha)=\max\limits_{i\in\alpha}f_i(\mathbf{c}_\alpha)=\langle\mathbf{e}_k,\mathbf{F}(\mathbf{c}_\alpha)\rangle.$$ This means that the Nash equilibrium $\mathbf{c}_\alpha$ is not the strictly Nash equilibrium for any $\alpha\subset \mathbf{I}_m$ with $|\alpha|\geq2$. As we already showed in \textbf{Part~$(ii)$} that the vertex $\mathbf{e}_k$ of the simplex is both stable and attracting for all $k\in\mathbf{I}_m$. Hence, a strictly Nash equilibrium is asymptotically stable.
			
			\item[Part $(iv)$] Due to \textbf{Part~$(ii)$}, an interior orbit starting from any initial point $\mathbf{x}\in \textup{int}\mathbb{S}^{m-1}$ converges to the center $\mathbf{c}_\beta$ of the face $\mathbb{S}_{\beta}$ where $\beta=\mathsf{\textbf{MaxInd}}(\mathbf{x})$ which is a Nash equilibrium. This completes the proof of Theorem \ref{FTEGT}.
		\end{description}
		
		\section{The Proof of Theorem \ref{HistircBehavior}}
		
		It is easy to check that two discrete-time replicator equations of a zero-sum game given by \eqref{REH1} and \eqref{REH2} resemble the similar dynamics. Hence, we discuss only the first discrete-time replicator equation given by \eqref{REH1}. We are aiming to show that 
		the discrete-time replicator equation $\mathcal{R}_{H}:\mathbb{S}^{2}\to\mathbb{S}^{2}$ given by \eqref{REH1} has \textit{uniformly historic behavior}. Recall (see \cite{Saburov2021d}) that the discrete-time replicator equation $\mathcal{R}_{H}:\mathbb{S}^{2}\to\mathbb{S}^{2}$ given by \eqref{REH1} is said to have \textit{uniformly historic behavior} if there exist at least two disjoint convex compact subsets $U_1,U_2\subset\mathbb{S}^{2}$ such that \textit{a uniformly historic initial set} $\mathcal{HIS}_{\mathcal{R}_{H}}(U_1,U_2)$ has positive Lebesgue measure where $\mathcal{HIS}_{\mathcal{R}_{H}}(U_1,U_2)$ is
		a set of all initial points $\mathbf{x}\in\mathbb{S}^{2}$ whose orbits have uniformly $(U_1,U_2)-$historic behavior, i.e., there exist $\lambda_i, \mu_i>0$ (depending on $\mathbf{x}$) for $i=1,2$ such that the trapping-escaping sequence 
		{
			\begin{multline*}
				\left\{\chi_{{U}_i}\left(\mathcal{R}_{H}^{(n)}(\mathbf{x})\right)\right\}_{n=0}^\infty=\left\{1^{p^{(i)}_{n}}0^{q^{(i)}_{n}}\right\}_{n=1}^\infty:=\\
				\left\{\underbrace{1,\cdots,1}_{p^{(i)}_1}, \ \underbrace{0,\cdots, 0}_{q^{(i)}_1}, \ \underbrace{1,\cdots,1}_{p^{(i)}_2}, \ \underbrace{0,\cdots,0}_{q^{(i)}_2}, \ \cdots \ \underbrace{1,\cdots,1}_{p^{(i)}_{n}}, \ \underbrace{0,\cdots,0}_{q^{(i)}_{n}},\cdots\right\}
			\end{multline*} 
		}
		of the orbit over the set $U_i$ has $(\lambda_i,\mu_i)-$gap for all $i=1,2$, i.e., 
		\begin{eqnarray*}
			\liminf\limits_{n\to\infty}\frac{p^{(i)}_{n+1}}{\sum\limits_{k=1}^{n}(p^{(i)}_k+q^{(i)}_k)}\geq\lambda_i \quad \textup{and} \quad \liminf\limits_{n\to\infty}\frac{q^{(i)}_{n+1}}{\sum\limits_{k=1}^{n}(p^{(i)}_k+q^{(i)}_k)+p^{(i)}_{n+1}}\geq\mu_i, \quad  i=1,2
		\end{eqnarray*} 
		where $\chi_U:U\to\mathbb{R}$ is a characteristic function of a set $U$, i.e., $\chi_U(\mathbf{x})=1$ if $\mathbf{x}\in U$ and $\chi_U(\mathbf{x})=0$ if $\mathbf{x}\not\in U$. 
		
		One of the main results of the paper \cite{Saburov2021d} (see Theorem B, \cite{Saburov2021d}) states that if the discrete-time replicator equation $\mathcal{R}_{H}:\mathbb{S}^{2}\to\mathbb{S}^{2}$ given by \eqref{REH1} has \textit{uniformly historic behavior} then for any $\mathbf{x}\in\mathcal{HIS}_{\mathcal{R}_{H}}(U_1,U_2)$, the $s^{th}-$order repeated time averages $\left\{\mathcal{A}^{(s)}_{n}\left(\mathbf{x}\right)\right\}_{n=1}^\infty$ of the discrete-time replicator equation $\mathcal{R}_{H}:\mathbb{S}^{2}\to\mathbb{S}^{2}$ do not converge for any $s\in\mathbb{N}$. That's why, in order to prove Theorem \ref{HistircBehavior},  it is enough to show that 
		the discrete-time replicator equation $\mathcal{R}_{H}:\mathbb{S}^{2}\to\mathbb{S}^{2}$ given by \eqref{REH1} has \textit{uniformly historic behavior}.
		
		\begin{description}
			\item[\textbf{Step-1}] Since $f_1(\mathbf{x}), f_2(\mathbf{x}), f_3(\mathbf{x})>0$ for any $\mathbf{x}\in\mathbb{S}^2$, it is easy to check that the rest (fixed) points of the discrete-time replicator equation \eqref{REH1} which lie on the boundary $\partial\mathbb{S}^2=\{\mathbf{x}\in\mathbb{S}^2: x_1x_2x_3=0\}$ of 2D simplex $\mathbb{S}^2$ are only three vertices $\mathbf{e}_1=(1,0,0)$, $\mathbf{e}_2=(0,1,0)$, and $\mathbf{e}_3=(0,0,1)$.  Moreover, since $f_1(\mathbf{c})=f_2(\mathbf{c})=f_3(\mathbf{c})$, the center $\mathbf{c}=(\frac{1}{3},\frac{1}{3},\frac{1}{3})$ of 2D simplex is an interior rest (fixed) point. We now show that there is no any other interior rest point except the center $\mathbf{c}=(\frac{1}{3},\frac{1}{3},\frac{1}{3})$ of the simplex $\mathbf{x}\in\mathbb{S}^2$. Let $\mathbf{x}\in\textup{int}\mathbb{S}^2$ be any interior rest (fixed) point. We must then have  $x_1f_1(\mathbf{x})=x_3f_2(\mathbf{x})$, $x_2f_2(\mathbf{x})=x_1f_3(\mathbf{x})$, and $x_3f_3(\mathbf{x})=x_2f_1(\mathbf{x})$. Consequently, we derive $$x_1f_1(\mathbf{x})+x_2f_2(\mathbf{x})+x_3f_3(\mathbf{x})=x_3f_2(\mathbf{x})+x_1f_3(\mathbf{x})+x_2f_1(\mathbf{x}).$$ 
			
			Since $f_i(\mathbf{x}) \left(\gtreqqless\right) f_j(\mathbf{x})$ if and only if $x_i\left(\gtreqqless\right)x_j$ for all $i,j\in\mathbf{I}_3$, it follows from the rearrangement inequality (see \cite{MOA}) that  $\mathbf{x}=\mathbf{c}$. 
			
			Three vertices of the simplex are the saddle points 
			which are attractive only along the boundaries of the simplex in the clockwise direction. Moreover, since a similar-order preserving mapping is continuously differentiable, the center of the simplex is a repelling rest (fixed) point with complex eigenvalues which induce a spiral behavior towards the boundary of the simplex. We are aiming to show that 
			the orbit of any interior initial point (except the center of the simplex)  
			alternatively spends exponentially increasing amount of times close to the consecutive vertices
			(by clockwise direction) of the simplex. 
			
			\item[\textbf{Step-2}]  We now  show that  $\xi(\mathbf{x}):=x_1x_2x_3$ is a decreasing Lyapunov function. Indeed, it follows from the arithmetic-geometric mean and rearrangement inequalities (see \cite{MOA}) that
			\begin{align*}
				\xi(\mathcal{R}_H(\mathbf{x}))\leq\xi(\mathbf{x})\left(1+\frac{\langle\pi(\mathbf{x}),\mathbf{F}(\mathbf{x})\rangle-\langle\mathbf{x},\mathbf{F}(\mathbf{x})\rangle}{3}\right)^{\frac{1}{3}}\leq\xi(\mathbf{x}).
			\end{align*}
			where $\pi(\mathbf{x}):=(x_2,x_3,x_1)$.
			Hence, $\left\{\xi\left(\mathcal{R}_H^{(n)}(\mathbf{x})\right)\right\}
			_{n=0}^\infty$
			is a decreasing sequence and converges to some limit $\lambda$. Since $\mathbf{x}\neq \mathbf{c}$, we have
			$$0\leq\lambda<\xi(\mathbf{x})<\xi(\mathbf{c}).$$
			If $\lambda>0$ then
			\begin{equation*}
				1=\lim\frac{\xi\big(\mathcal{R}_H^{(n+1)}(\mathbf{x})\big)}
				{\xi\big(\mathcal{R}_H^{(n)}(\mathbf{x})\big)}
				=\lim\zeta\big(\mathcal{R}_H^{(n)}(\mathbf{x})\big)\leq 1,
			\end{equation*}
			where
			\begin{multline*}
				\zeta(\mathbf{x})=(1+x_2f_1(\mathbf{x})-x_3f_3(\mathbf{x}))
				(1+x_3f_2(\mathbf{x})-x_1f_1(\mathbf{x}))(1+x_1f_3(\mathbf{x})-x_2f_2(\mathbf{x})).
			\end{multline*}
			Consequently, for any $\mathbf{x}^{*}\in\omega(\mathbf{x})$ (where $\omega(\mathbf{x})$ is an omega-limit set of an orbit starting from an initial point~$\mathbf{x}$), 
			one has 
			$1=\zeta(\mathbf{x}^{*})=
			\max\limits_{\mathbf{y}\in\mathbb{S}^2}\zeta(\mathbf{y})=1$, 
			i.e,   $\mathbf{x}^{*}=\mathbf{c}$. 
			However, this contradicts to
			$\xi(\mathbf{x}^{*})=\lambda<\xi(\mathbf{c})$. It shows that $\lambda=~0$, i.e., 
			$\omega(\mathbf{x})\subset\partial\mathbb{S}^2$. 
			
			Obviously, we have
			$\omega(\mathbf{x})\neq\{\mathbf{e}_1,\mathbf{e}_2,\mathbf{e}_3\}$. 
			Since $\mathcal{R}_H(\omega(\mathbf{x}))=\omega(\mathbf{x})$, the orbit  
			$\{\mathcal{R}_H^{(n)}(\mathbf{x}^{*})\}_{n=0}^\infty
			\subset \omega(\mathbf{x})$  
			is infinite for any
			$\mathbf{x}^{*}\in
			\omega(\mathbf{x})\setminus\{\mathbf{e}_1,\mathbf{e}_2,\mathbf{e}_3\}$. It means that an omega-limit set $\omega(\mathbf{x})$ of an orbit  
			starting from an initial point 
			$\mathbf{x}\in {\textup{int}}\mathbb{S}^{2}\setminus \{\mathbf{c}\}$ 
			is infinite and lies on the boundary 
			$\partial\mathbb{S}^{2}$ of the simplex $\mathbb{S}^{2}$. 
			
			\item[\textbf{Step-3}]
			An orbit starting from any interior initial point eventually 
			moves along the periodic itinerary 
			$$G_1\hookrightarrow G_2\hookrightarrow G_3\hookrightarrow
			G_4\hookrightarrow G_5\hookrightarrow G_6\hookrightarrow G_1$$ where $G_i\hookrightarrow G_j$ 
			stands for $\mathcal{R}_H(G_i)\subset G_i\cup G_j$ and the set $G_i$ is defined as follows for $i=1,2\cdots,6$ 
			{\begin{align*}
					&G_1=\Big\{\mathbf{x}\in\mathbb{S}^2: x_1\geq x_2\geq
					x_3\Big\}, \quad
					G_2=\Big\{\mathbf{x}\in\mathbb{S}^2: x_1\geq x_3\geq
					x_2\Big\},\\
					&G_3=\Big\{\mathbf{x}\in\mathbb{S}^2: x_3\geq x_1\geq 
					x_2\Big\}, \quad 
					G_4=\Big\{\mathbf{x}\in\mathbb{S}^2: x_3\geq x_2\geq
					x_1\Big\},\\
					&G_5=\Big\{\mathbf{x}\in\mathbb{S}^2: x_2 \geq x_3\geq
					x_1\Big\}, \quad
					G_6=\Big\{\mathbf{x}\in\mathbb{S}^2: x_2\geq x_1\geq
					x_3\Big\}.
			\end{align*}}
			
			Let us show that $G_1\hookrightarrow G_2$.  The rest is similar to this case. Let 
			$\partial\mathbb{S}^2_\varepsilon=\{\mathbf{x}\in \mathbb{S}^2: 
			{\textup{dist}}(\mathbf{x},\partial\mathbb{S}^2)<\varepsilon\}$. 
			Since $\omega(\mathbf{x})\subset \partial\mathbb{S}^2$, 
			one has  
			$\mathcal{R}_H^{(n)}(\mathbf{x})\in\partial\mathbb{S}^2_\varepsilon$
			for sufficiently small $\varepsilon>0$ and $n\geq n_0$. 
			We show 
			$\mathcal{R}_H(G_1\cap\partial\mathbb{S}^2_\varepsilon)
			\subset G_1\cup G_2$.
			
			Let $\mathbf{x}=(x_1,x_2,x_3)\in G_1\cap\partial\mathbb{S}^2_\varepsilon$ and $\mathbf{x}^{(1)}=(x_1^{(1)},x_2^{(1)},x_3^{(1)})=\mathcal{R}_H(\mathbf{x})$.
			Since
			$\mathbf{x}\in G_1$, we have  $x_1\geq x_2\geq
			x_3$. Moreover, since $f_i(\mathbf{x}) \left(\gtreqqless\right) f_j(\mathbf{x})$ if and only if $x_i\left(\gtreqqless\right)x_j$ for all $i,j\in\mathbf{I}_3$,
			we obtain
			\begin{equation*}
				x_2f_1(\mathbf{x})-x_3f_3(\mathbf{x})\geq0 \quad \textup{and} \quad x_3f_2(\mathbf{x})-x_1f_1(\mathbf{x})\leq0.
			\end{equation*}
			Consequently, we have ${{x}^{(1)}_1}\geq {x_1}\geq 
			{x}_2\geq x^{(1)}_2.$
			Obviously, for
			$\mathbf{x}\in G_1\cap\partial\mathbb{S}^2_\varepsilon$,
			one has $x_3<\varepsilon$ for sufficiently small~$\varepsilon$. We then obtain  
			${x^{(1)}_3}\leq 2x_3\leq 
			2\varepsilon\leq
			{x_1}\leq{x^{(1)}_1}.$
			Therefore, we obtain
			$\mathcal{R}_H(G_1\cap\partial\mathbb{S}^2_\varepsilon)
			\subset G_1\cup G_2$.
			
			\item[\textbf{Step-4}] We choose a neighborhood $\mathbb{U}_0$ 
			of the repelling rest point $\mathbf{c}$ 
			such that $\mathbb{U}_0\subset \textup{int}\mathbb{S}^2$ and 
			$\mathbb{U}_1=(G_1\cup G_2)\setminus \mathbb{U}_0$, \ 
			$\mathbb{U}_2=(G_3\cup G_4)\setminus \mathbb{U}_0$, \
			$\mathbb{U}_3=(G_5\cup G_6)\setminus \mathbb{U}_0$ 
			are convex compact sets which satisfy 
			$\mathbb{U}_1\cap \mathbb{U}_2\cap \mathbb{U}_3=\emptyset$
			
			Let $\mathbb{U}$ be one of the sets 
			$\mathbb{U}_1$, $\mathbb{U}_2$, and $\mathbb{U}_3$ defined above and let $\mathbf{x}\notin \mathbb{U}\cup \mathbb{U}_0$,
			$\mathcal{R}_H^{(k)}(\mathbf{x})\in \mathbb{U}$ 
			for all $1\leq k\leq n$, 
			and $\mathcal{R}_H^{(n+1)}(\mathbf{x})\notin \mathbb{U}$.
			Then we have 
			$n> A\log\frac{B}{\xi(\mathcal{R}_H(\mathbf{x}))}$  
			where $A$ and $B$ are some positive constants.
			
			Indeed, let $\mathbf{x}=(x_1,x_2,x_3)$ and
			$\mathbf{x}^{(k)}=(x^{(k)}_1,x_2^{(k)},x_3^{(k)})=\mathcal{R}_H^{(k)}(\mathbf{x})$. We also assume $\mathbb{U}:=\mathbb{U}_1$. 
			Since $\mathbf{x}\notin \mathbb{U}_1$,
			$\mathcal{R}_H^{(k)}(\mathbf{x})\in \mathbb{U}_1$ 
			for all $1\leq k\leq n$, 
			and $\mathcal{R}_H^{(n+1)}(\mathbf{x})\notin \mathbb{U}_1$,
			we then obtain $\mathbf{x}\in G_6$,
			$\mathcal{R}_H(\mathbf{x})\in G_1$, and
			$\mathcal{R}_H^{(n+1)}(\mathbf{x})\in G_3$.
			Therefore, $x_2\geq\frac{1}{3}$, $x^{(1)}_1\geq\frac{1}{3}$, and 
			$x^{(n+1)}_3\geq \frac{1}{3}$. 
			Hence, we have that
			$$x_2^{(1)}=x_2\left(1+x_3f_2(\mathbf{x})-x_1f_1(\mathbf{x})\right)
			\geq x_2^2\geq\frac{1}{9}$$ and
			\begin{align*}
				\frac{1}{81\xi(\mathcal{R}_H(\mathbf{x}))}
				\leq
				\frac{x_1^{(1)}x_2^{(1)}}
				{3\xi(\mathcal{R}_H(\mathbf{x}))}=
				\frac{1}{3x_3^{(1)}}\leq
				\frac{x_3^{(n+1)}}{x_3^{(1)}}=
				\prod\limits_{k=1}^{n}\frac{x_3^{(k+1)}}{x_3^{(k)}}\leq 2^n,
			\end{align*}
			which yields 
			$n> A\log\frac{B}{\xi(\mathcal{R}_H(\mathbf{x}))}$ 
			for some constants $A>0$ and $B>0$. 
			
			\item[\textbf{Step-5}] Let $\mathbb{U}$ be one of the sets 
			$\mathbb{U}_1$, $\mathbb{U}_2$, and $\mathbb{U}_3$ given above and 
			let $(n_i,m_i)_{i=1}^\infty$ be a sequence of natural numbers 
			such that $\mathcal{R}_H^{(n_i)}(\mathbf{x})\notin \mathbb{U}$,
			$\mathcal{R}_H^{(n_i+k)}(\mathbf{x})\in \mathbb{U}$ 
			for $1\leq k\leq m_i$, and
			$\mathcal{R}_H^{(n_i+m_i+1)}(\mathbf{x})\notin \mathbb{U}$.

			It follows from \textbf{Step-2} and \textbf{Step-3} that $n_i\to \infty$.
			Let $\rho:=\max\limits_
			{\mathbf{x}\in\mathbb{S}^2\setminus \mathbb{U}_0}\zeta(\mathbf{x})<1$ (see \textbf{Step-2}). 
			Then, we obtain that 
			$$\xi(\mathcal{R}_H(\mathbf{x}))=
			\xi(\mathbf{x})\zeta(\mathbf{x})<\rho\xi(\mathbf{x})$$ 
			for any $\mathbf{x}\in\mathbb{S}^2\setminus \mathbb{U}_0$.  
			It follows from \textbf{Step-4} that
			$$m_i>A\log\frac{B}{\xi(\mathcal{R}_H^{(n_i)}(\mathbf{x}))}>
			A\log\frac{B}{\rho^{n_i}\xi(\mathbf{c})}>Cn_i$$ for some constant $C>0$.
			
			\item[\textbf{Step-6}] It follows from \textbf{Step-3} that the forward orbit eventually moves along the periodic itinerary
			\begin{multline*}
				\underbrace{{\mathbb{U}_1}\ {^{\curvearrowright}}\ {\mathbb{U}_1}}_{p_1}\rightarrow\underbrace{{\mathbb{U}_2}\ {^{\curvearrowright}}\ {\mathbb{U}_2}}_{q_1}\rightarrow\underbrace{{\mathbb{U}_3}\ {^{\curvearrowright}}\ {\mathbb{U}_3}}_{r_1}\rightarrow\\
				\underbrace{{\mathbb{U}_1}\ {^{\curvearrowright}}\ {\mathbb{U}_1}}_{p_2}\rightarrow\underbrace{{\mathbb{U}_2}\ {^{\curvearrowright}}\ {\mathbb{U}_2}}_{q_2}\rightarrow\underbrace{{\mathbb{U}_3}\ {^{\curvearrowright}}\ {\mathbb{U}_3}}_{r_2}\rightarrow\\
				\cdots\rightarrow\underbrace{{\mathbb{U}_1}\ {^{\curvearrowright}}\ {\mathbb{U}_1}}_{p_n}\rightarrow\underbrace{{\mathbb{U}_2}\ {^{\curvearrowright}}\ {\mathbb{U}_2}}_{q_n}\rightarrow\underbrace{{\mathbb{U}_3}\ {^{\curvearrowright}}\ {\mathbb{U}_3}}_{r_n}\rightarrow\cdots
			\end{multline*}
			where $\{p_n\}$, $\{q_n\}$, and $\{r_n\}$ are sequences of trapping times of the forward orbit over the sets  $\mathbb{U}_1$, $\mathbb{U}_2$, and $\mathbb{U}_3$, respectively.
			It follows from \textbf{Step-5} that there exists a constant $C>0$ such that
			\begin{multline*}
				p_{n+1}>C\sum_{i=1}^{n}(p_i+q_i+r_i),  \\ q_{n+1}>C\left(\sum_{i=1}^{n}(p_i+q_i+r_i)+p_{n+1}\right), \\ r_{n+1}>C\left(\sum_{i=1}^{n}(p_i+q_i+r_i)+p_{n+1}+q_{n+1}\right).
			\end{multline*}
			Hence, the trapping-escaping sequence $\left\{\chi_{\mathbb{U}_i}(\mathcal{R}_H^{(n)}\mathbf{x})\right\}_{n=0}^\infty$ has $(C;C(C+2))-$ gap (see \cite{Saburov2021d}) for $i=1,2,3$ where $\chi_\mathbb{U}(\cdot)$ is a characteristic function of a set $\mathbb{U}$, i.e., $\chi_\mathbb{U}(\mathbf{x})=1$ if $\mathbf{x}\in \mathbb{U}$ and $\chi_\mathbb{U}(\mathbf{x})=0$ if $\mathbf{x}\not\in \mathbb{U}$. Consequently, the discrete-time replicator equation of a zero-sum game given by \eqref{REH1} has \textit{uniformly historic behavior} which eventually causes the divergence of higher-order repeated time averages $\left\{\mathcal{A}^{(s)}_{n}\left(\mathbf{x}\right)\right\}_{n=1}^\infty$ for any $s\in\mathbb{N}$ where $\mathcal{A}^{(s)}_{n}\left(\mathbf{x}\right):=\frac{1}{n}\sum_{k=1}^{n}\mathcal{A}^{(s-1)}_{k}\left(\mathbf{x}\right)$ and $\mathcal{A}^{(1)}_{n}\left(\mathbf{x}\right):=\frac{1}{n}\sum_{k=0}^{n-1}\mathcal{R}_{H}^{(k)}(\mathbf{x})$ (see Theorem B, \cite{Saburov2021d}). It completes the proof of Theorem \ref{HistircBehavior}.

		\end{description}

		\bibliographystyle{amsplain}
		\bibliography{paper}

\end{document}